%Header

\magnification1200
\baselineskip15pt
\let\otheta=\theta
%\input makros.mak
%---Anfang makros.mak
\newread\AUX\immediate\openin\AUX=\jobname.aux
\newcount\relFnno
\def\ref#1{\expandafter\edef\csname#1\endcsname}
\ifeof\AUX\immediate\write16{\jobname.aux gibt es nicht!}\else
\input \jobname.aux
\fi\immediate\closein\AUX

% Wirkungsweise:

% ...bearbeiteter Text...
% \ignore                   <--- muss am Ende einer Zeile stehen
% ...ignorierter Text...    <--- keine Zeile darf mit ';' beginnen
% ;;                        <--- muss allein auf einer Zeile stehen
% ...bearbeiteter Text...

\def\ignore{\bgroup
\catcode`\;=0\catcode`\^^I=14\catcode`\^^J=14\catcode`\^^M=14
\catcode`\ =14\catcode`\!=14\catcode`\"=14\catcode`\#=14\catcode`\$=14
\catcode`\&=14\catcode`\'=14\catcode`\(=14\catcode`\)=14\catcode`\*=14
\catcode`+=14\catcode`\,=14\catcode`\-=14\catcode`\.=14\catcode`\/=14
\catcode`\0=14\catcode`\1=14\catcode`\2=14\catcode`\3=14\catcode`\4=14
\catcode`\5=14\catcode`\6=14\catcode`\7=14\catcode`\8=14\catcode`\9=14
\catcode`\:=14\catcode`\<=14\catcode`\==14\catcode`\>=14\catcode`\?=14
\catcode`\@=14\catcode`\A=14\catcode`\B=14\catcode`\C=14\catcode`\D=14
\catcode`\E=14\catcode`\F=14\catcode`\G=14\catcode`\H=14\catcode`\I=14
\catcode`\J=14\catcode`\K=14\catcode`\L=14\catcode`\M=14\catcode`\N=14
\catcode`\O=14\catcode`\P=14\catcode`\Q=14\catcode`\R=14\catcode`\S=14
\catcode`\T=14\catcode`\U=14\catcode`\V=14\catcode`\W=14\catcode`\X=14
\catcode`\Y=14\catcode`\Z=14\catcode`\[=14\catcode`\\=14\catcode`\]=14
\catcode`\^=14\catcode`\_=14\catcode`\`=14\catcode`\a=14\catcode`\b=14
\catcode`\c=14\catcode`\d=14\catcode`\e=14\catcode`\f=14\catcode`\g=14
\catcode`\h=14\catcode`\i=14\catcode`\j=14\catcode`\k=14\catcode`\l=14
\catcode`\m=14\catcode`\n=14\catcode`\o=14\catcode`\p=14\catcode`\q=14
\catcode`\r=14\catcode`\s=14\catcode`\t=14\catcode`\u=14\catcode`\v=14
\catcode`\w=14\catcode`\x=14\catcode`\y=14\catcode`\z=14\catcode`\{=14
\catcode`\|=14\catcode`\}=14\catcode`\~=14\catcode`\^^?=14
\Ignoriere}
\def\Ignoriere#1\;{\egroup}

\newcount\itemcount
\def\resetitem{\global\itemcount0}\resetitem
\newcount\Itemcount
\Itemcount0
\newcount\maxItemcount
\maxItemcount=0

\def\FILTER\fam\itfam\tenit#1){#1}

\def\Item#1{\global\advance\itemcount1
\edef\TEXT{{\it\romannumeral\itemcount)}}%
\ifx?#1?\relax\else
\ifnum#1>\maxItemcount\global\maxItemcount=#1\fi
\expandafter\ifx\csname I#1\endcsname\relax\else
\edef\testA{\csname I#1\endcsname}
\expandafter\expandafter\def\expandafter\testB\testA
\edef\testC{\expandafter\FILTER\testB}
\edef\testD{\csname0\testC0\endcsname}\fi
\edef\testE{\csname0\romannumeral\itemcount0\endcsname}
\ifx\testD\testE\relax\else
\immediate\write16{I#1 hat sich geaendert!}\fi
\expandwrite\AUX{\neverexpand\ref{I#1}{\TEXT}}\fi
\item{\TEXT}}

\def\today{\number\day.~\ifcase\month\or
  Januar\or Februar\or M{\"a}rz\or April\or Mai\or Juni\or
  Juli\or August\or September\or Oktober\or November\or Dezember\fi
  \space\number\year}
\font\sevenex=cmex7
%\font\sevenex=cmex10 scaled 700
\scriptfont3=\sevenex
%\font\fiveex=cmex7 scaled 714
\font\fiveex=cmex10 scaled 500
\scriptscriptfont3=\fiveex
\def\A{{\bf A}}
\def\G{{\bf G}}
\def\P{{\bf P}}

\def\phi{\varphi}
\def\epsilon{\varepsilon}
\def\theta{\vartheta}
%\font\lams=lams3
%\def\uinto{\lower1.7pt\hbox{%
%\vbox{\offinterlineskip
%\hbox{\lams\char"7A}%
%\hbox{\vbox to 7.5pt{\leaders\vrule width0.2pt\vfill}%
%\kern-.3pt\hbox{\lams\char"76}}}}}
\def\uauf{\lower1.7pt\hbox to 3pt{%
\vbox{\offinterlineskip
\hbox{\vbox to 8.5pt{\leaders\vrule width0.2pt\vfill}%
\kern-.3pt\hbox{\lams\char"76}\kern-0.3pt%
$\raise1pt\hbox{\lams\char"76}$}}\hfil}}

\def\title#1{\par
{\baselineskip1.5\baselineskip\rightskip0pt plus 5truecm
\leavevmode\vskip0truecm\noindent\font\BF=cmbx10 scaled \magstep2\BF #1\par}
\vskip1truecm
\leftline{\font\CSC=cmcsc10
{\CSC Friedrich Knop}\footnote{This work originates from a stay at the University of Strasbourg in 1996 and was finished during a stay at the University of Freiburg in 2003. The author thanks both institutions for their hospitality.}}
\leftline{Department of Mathematics, Rutgers University, New Brunswick NJ
08903, USA}
\leftline{knop@math.rutgers.edu}
\vskip1truecm
\par}

%%%%%%%%%%%%%%%%%%
% Makros f"ur Querverweise:
\def\cite#1{\expandafter\ifx\csname#1\endcsname\relax
{\bf?}\immediate\write16{#1 ist nicht definiert!}\else\csname#1\endcsname\fi}
\def\expandwrite#1#2{\edef\next{\write#1{#2}}\next}
\def\neverexpand{\noexpand\noexpand\noexpand}
\def\strip#1\ {}
\def\ncite#1{\expandafter\ifx\csname#1\endcsname\relax
{\bf?}\immediate\write16{#1 ist nicht definiert!}\else
\expandafter\expandafter\expandafter\strip\csname#1\endcsname\fi}
\newwrite\AUX
\immediate\openout\AUX=\jobname.aux
%%%%%%%%%%%%%%%%%%%%%
% Definition von \eightpoint:
\font\eightrm=cmr8\font\sixrm=cmr6
\font\eighti=cmmi8
\font\eightit=cmti8
\font\eightbf=cmbx8
\font\eightcsc=cmcsc10 scaled 833
\def\eightpoint{%
\textfont0=\eightrm\scriptfont0=\sixrm\def\rm{\fam0\eightrm}%
\textfont1=\eighti
\textfont\bffam=\eightbf\def\bf{\fam\bffam\eightbf}%
\textfont\itfam=\eightit\def\it{\fam\itfam\eightit}%
\def\csc{\eightcsc}%
\setbox\strutbox=\hbox{\vrule height7pt depth2pt width0pt}%
\normalbaselineskip=0,8\normalbaselineskip\normalbaselines\rm}
%%%%%%%%%%%%%%%%%%
% Fussnotenmakros
\newcount\absFnno\absFnno1
\write\AUX{\relFnno1}
\newif\ifMARKE\MARKEtrue
{\catcode`\@=11
\gdef\footnote{\ifMARKE\edef\@sf{\spacefactor\the\spacefactor}\/%
$^{\cite{Fn\the\absFnno}}$\@sf\fi
\MARKEtrue
\insert\footins\bgroup\eightpoint
\interlinepenalty100\let\par=\endgraf
\leftskip=0pt\rightskip=0pt
\splittopskip=10pt plus 1pt minus 1pt \floatingpenalty=20000\smallskip
\item{$^{\cite{Fn\the\absFnno}}$}%
\expandwrite\AUX{\neverexpand\ref{Fn\the\absFnno}{\neverexpand\the\relFnno}}%
\global\advance\absFnno1\write\AUX{\advance\relFnno1}%
\bgroup\strut\aftergroup\@foot\let\next}}
\skip\footins=12pt plus 2pt minus 4pt
\dimen\footins=30pc
\output={\plainoutput\immediate\write\AUX{\relFnno1}}
%%%%%%%%%%%%%%%%%%%%%
\newcount\Abschnitt\Abschnitt0
\def\beginsection#1. #2 \par{\advance\Abschnitt1%
\vskip0pt plus.10\vsize\penalty-250
\vskip0pt plus-.10\vsize\bigskip\vskip\parskip
\edef\TEST{\number\Abschnitt}
\expandafter\ifx\csname#1\endcsname\TEST\relax\else
\immediate\write16{#1 hat sich geaendert!}\fi
\expandwrite\AUX{\neverexpand\ref{#1}{\TEST}}
\leftline{\marginnote{#1}\bf\number\Abschnitt. \ignorespaces#2}%
\nobreak\smallskip\noindent\SATZ1\GNo0}
%%%%%%%%%%%%%%%%%%
\def\Proof:{\par\noindent{\it Proof:}}
\def\Remark:{\ifdim\lastskip<\medskipamount\removelastskip\medskip\fi
\noindent{\bf Remark:}}
\def\Remarks:{\ifdim\lastskip<\medskipamount\removelastskip\medskip\fi
\noindent{\bf Remarks:}}
\def\Definition:{\ifdim\lastskip<\medskipamount\removelastskip\medskip\fi
\noindent{\bf Definition:}}
\def\Example:{\ifdim\lastskip<\medskipamount\removelastskip\medskip\fi
\noindent{\bf Example:}}
\def\Examples:{\ifdim\lastskip<\medskipamount\removelastskip\medskip\fi
\noindent{\bf Examples:}}
%%%%%%%%%%%%%%%%
\newif\ifmarginalnotes\marginalnotesfalse
\newif\ifmarginalwarnings\marginalwarningstrue

\def\marginnote#1{\ifmarginalnotes\hbox to 0pt{\eightpoint\hss #1\ }\fi}

\def\strutdepth{\dp\strutbox}
\def\Randbem#1#2{\ifmarginalwarnings
{#1}\strut
\setbox0=\vtop{\eightpoint
\rightskip=0pt plus 6mm\hfuzz=3pt\hsize=16mm\noindent\leavevmode#2}%
\vadjust{\kern-\strutdepth
\vtop to \strutdepth{\kern-\ht0
\hbox to \hsize{\kern-16mm\kern-6pt\box0\kern6pt\hfill}\vss}}\fi}

\def\Zitat!{\Randbem{\bf?}{\bf Zitat}}

\newcount\SATZ\SATZ1
\def\proclaim #1. #2\par{\ifdim\lastskip<\medskipamount\removelastskip
\medskip\fi
\noindent{\bf#1.\ }{\it#2}\Par
\ifdim\lastskip<\medskipamount\removelastskip\goodbreak\medskip\fi}
\def\Aussage#1{\expandafter\def\csname#1\endcsname##1.{\resetitem
\ifx?##1?\relax\else
\edef\TEST{#1\penalty10000\ \number\Abschnitt.\number\SATZ}
\expandafter\ifx\csname##1\endcsname\TEST\relax\else
\immediate\write16{##1 hat sich geaendert!}\fi
\expandwrite\AUX{\neverexpand\ref{##1}{\TEST}}\fi
\proclaim {\marginnote{##1}\number\Abschnitt.\number\SATZ. #1\global\advance\SATZ1}.}}
\Aussage{Theorem}
\Aussage{Proposition}
\Aussage{Corollary}
\Aussage{Lemma}
%%%%%%%%%%%%%%%%
\font\la=lasy10
\def\strich{\hbox{$\vcenter{\hbox
to 1pt{\leaders\hrule height -0,2pt depth 0,6pt\hfil}}$}}
\def\dashedrightarrow{\hbox{%
\hbox to 0,5cm{\leaders\hbox to 2pt{\hfil\strich\hfil}\hfil}%
\kern-2pt\hbox{\la\char\string"29}}}

\def\Bindestrich{\penalty10000-\hskip0pt}
\let\_=\Bindestrich
\def\.{{\sfcode`.=1000.}}
%%%%%%%%%%%%%%%%%%%%%%%%%%%%%%%%%%%%

\def\Par{\par}
\def\:={\mathrel{\raise0,9pt\hbox{.}\kern-2,77779pt
\raise3pt\hbox{.}\kern-2,5pt=}}
\def\=:{\mathrel{=\kern-2,5pt\raise0,9pt\hbox{.}\kern-2,77779pt
\raise3pt\hbox{.}}} 

\def\pfeil{\rightarrow}

\def\Pfeil{\longrightarrow}
\def\pf#1{\buildrel#1\over\rightarrow}
\def\Pf#1{\buildrel#1\over\longrightarrow}

\def\Ugleich{\hbox{$\cup$\kern.5pt\vrule depth -0.5pt}}
\def\|#1|{\mathop{\rm#1}\nolimits}
\def\<{\langle}
\def\>{\rangle}
\let\Times=\times
\def\times{\mathop{\Times}}
\let\Otimes=\otimes
\def\otimes{\mathop{\Otimes}}
%%%%%%%%%%%%%%%%%%%%%%%%%%%%%%%%%
%Laden von Fonts:
\catcode`\@=11
\def\hex#1{\ifcase#1 0\or1\or2\or3\or4\or5\or6\or7\or8\or9\or A\or B\or
C\or D\or E\or F\else\message{Warnung: Setze hex#1=0}0\fi}
\def\fontdef#1:#2,#3,#4.{%
\alloc@8\fam\chardef\sixt@@n\FAM
\ifx!#2!\else\expandafter\font\csname text#1\endcsname=#2
\textfont\the\FAM=\csname text#1\endcsname\fi
\ifx!#3!\else\expandafter\font\csname script#1\endcsname=#3
\scriptfont\the\FAM=\csname script#1\endcsname\fi
\ifx!#4!\else\expandafter\font\csname scriptscript#1\endcsname=#4
\scriptscriptfont\the\FAM=\csname scriptscript#1\endcsname\fi
\expandafter\edef\csname #1\endcsname{\fam\the\FAM\csname text#1\endcsname}
\expandafter\edef\csname hex#1fam\endcsname{\hex\FAM}}
\catcode`\@=12 

%%%%%%%%%%%%%%%%%%%%%%%%%%%%%%%%%
\fontdef Ss:cmss10,,.
\fontdef Fr:eufm10,eufm7,eufm5.

                        %Hier aufpassen!!!

%
\fontdef bbb:msbm10,msbm7,msbm5.
\fontdef mbf:cmmib10,cmmib7,.
\fontdef msa:msam10,msam7,msam5.

\def\NN{{\bbb N}}

\def\ZZ{{\bbb Z}}

\def\cH{{\cal H}}
\def\cL{{\cal L}}
\def\cM{{\cal M}}

\mathchardef\leer=\string"0\hexbbbfam3F
\mathchardef\subsetneq=\string"3\hexbbbfam24
\mathchardef\semidir=\string"2\hexbbbfam6E
\mathchardef\dirsemi=\string"2\hexbbbfam6F
\mathchardef\haken=\string"2\hexmsafam78
\mathchardef\auf=\string"3\hexmsafam10
\let\OL=\overline
\def\overline#1{{\hskip1pt\OL{\hskip-1pt#1\hskip-.3pt}\hskip.3pt}}

%<--                    Aufpassen  

\def\vq{{\overline{v}}}
\def\wq{{\overline{w}}}
\def\xq{{\overline{x}}}

%
%%%%%%%%%%%%
\newdimen\Parindent
\Parindent=\parindent

%%%%%%%%%%%%
% Displayroutine

\abovedisplayskip 9.0pt plus 3.0pt minus 3.0pt
\belowdisplayskip 9.0pt plus 3.0pt minus 3.0pt
\newdimen\Grenze\Grenze2\Parindent\advance\Grenze1em
\newdimen\Breite
\newbox\DpBox
\def\NewDisplay#1
#2$${\Breite\hsize\advance\Breite-\hangindent
\setbox\DpBox=\hbox{\hskip2\Parindent$\displaystyle{\eqno{#1}#2}$}%
\ifnum\predisplaysize<\Grenze\abovedisplayskip\abovedisplayshortskip
\belowdisplayskip\belowdisplayshortskip\fi
\global\futurelet\nexttok\WEITER}
\def\WEITER{\ifx\nexttok\qed\expandafter\leftQEDdisplay
\else\leftdisplay\fi}
\def\leftdisplay{\hskip-\hangindent\leftline{\box\DpBox}$$}
\def\leftQEDdisplay{\hskip-\hangindent
\line{\copy\DpBox\hfill\lower\dp\DpBox\copy\QEDbox}%
\belowdisplayskip0pt$$\bigskip\let\nexttok=}
\everydisplay{\NewDisplay}
%%%%%%%%%%%
\newcount\GNo\GNo=0
\newcount\maxEqNo\maxEqNo=0
\def\eqno#1{%
\global\advance\GNo1
\edef\FTEST{$(\number\Abschnitt.\number\GNo)$}
\ifx?#1?\relax\else
\ifnum#1>\maxEqNo\global\maxEqNo=#1\fi%
\expandafter\ifx\csname E#1\endcsname\FTEST\relax\else
\immediate\write16{E#1 hat sich geaendert!}\fi
\expandwrite\AUX{\neverexpand\ref{E#1}{\FTEST}}\fi
\llap{\hbox to 40pt{\marginnote{#1}\FTEST\hfill}}}

\catcode`@=11
\def\eqalignno#1{\null\vcenter{\openup\jot\m@th\ialign{\eqno{##}\hfil
&\strut\hfil$\displaystyle{##}$&$\displaystyle{{}##}$\hfil\crcr#1\crcr}}\,}
\catcode`@=12

%%%%%%%%%%%%
\newbox\QEDbox
\newbox\nichts\setbox\nichts=\vbox{}\wd\nichts=2mm\ht\nichts=2mm
\setbox\QEDbox=\hbox{\vrule\vbox{\hrule\copy\nichts\hrule}\vrule}
\def\qed{\leavevmode\unskip\hfil\null\nobreak\hfill\copy\QEDbox\medbreak}
%%%%%%%%%%%%%%
\newdimen\HIindent
\newbox\HIbox
\def\setHI#1{\setbox\HIbox=\hbox{#1}\HIindent=\wd\HIbox}
\def\HI#1{\par\hangindent\HIindent\hangafter=0\noindent\leavevmode
\llap{\hbox to\HIindent{#1\hfil}}\ignorespaces}
%%%%%%%%%%%%%%

\newdimen\maxSpalbr
\newdimen\altSpalbr
\newcount\Zaehler

%{\catcode`/=\active
%\gdef\SlashOn{\catcode`/=\active\def/{X\string/\hskip0pt Y}}
%}

\newif\ifxxx

{\catcode`/=\active

\gdef\beginrefs{%
\xxxfalse
\catcode`/=\active
\def/{\string/\ifxxx\hskip0pt\fi}
\def\TText##1{{\xxxtrue\tt##1}}
\expandafter\ifx\csname Spaltenbreite\endcsname\relax
\def\Spaltenbreite{1cm}\immediate\write16{Spaltenbreite undefiniert!}\fi
\expandafter\altSpalbr\Spaltenbreite
\maxSpalbr0pt
\gdef\alt{}
\def\\##1\relax{%
\gdef\neu{##1}\ifx\alt\neu\global\advance\Zaehler1\else
\xdef\alt{\neu}\global\Zaehler=1\fi\xdef\SigText{##1\the\Zaehler}}
\def\L|Abk:##1|Sig:##2|Au:##3|Tit:##4|Zs:##5|Bd:##6|S:##7|J:##8|xxx:##9||{%
\def\SigText{##2}\global\setbox0=\hbox{##2\relax}
\edef\TEST{[\SigText]}
\expandafter\ifx\csname##1\endcsname\TEST\relax\else
\immediate\write16{##1 hat sich geaendert!}\fi
\expandwrite\AUX{\neverexpand\ref{##1}{\TEST}}
\setHI{[\SigText]\ }
\ifnum\HIindent>\maxSpalbr\maxSpalbr\HIindent\fi
\ifnum\HIindent<\altSpalbr\HIindent\altSpalbr\fi
\HI{\marginnote{##1}[\SigText]}
\ifx-##3\relax\else{##3}: \fi
\ifx-##4\relax\else{##4}{\sfcode`.=3000.} \fi
\ifx-##5\relax\else{\it ##5\/} \fi
\ifx-##6\relax\else{\bf ##6} \fi
\ifx-##8\relax\else({##8})\fi
\ifx-##7\relax\else, {##7}\fi
\ifx-##9\relax\else, \TText{##9}\fi\Par}
\def\B|Abk:##1|Sig:##2|Au:##3|Tit:##4|Reihe:##5|Verlag:##6|Ort:##7|J:##8|xxx:##9||{%
\def\SigText{##2}\global\setbox0=\hbox{##2\relax}
\edef\TEST{[\SigText]}
\expandafter\ifx\csname##1\endcsname\TEST\relax\else
\immediate\write16{##1 hat sich geaendert!}\fi
\expandwrite\AUX{\neverexpand\ref{##1}{\TEST}}
\setHI{[\SigText]\ }
\ifnum\HIindent>\maxSpalbr\maxSpalbr\HIindent\fi
\ifnum\HIindent<\altSpalbr\HIindent\altSpalbr\fi
\HI{\marginnote{##1}[\SigText]}
\ifx-##3\relax\else{##3}: \fi
\ifx-##4\relax\else{##4}{\sfcode`.=3000.} \fi
\ifx-##5\relax\else{(##5)} \fi
\ifx-##7\relax\else{##7:} \fi
\ifx-##6\relax\else{##6}\fi
\ifx-##8\relax\else{ ##8}\fi
\ifx-##9\relax\else, \TText{##9}\fi\Par}
\def\Pr|Abk:##1|Sig:##2|Au:##3|Artikel:##4|Titel:##5|Hgr:##6|Reihe:{%
\def\SigText{##2}\global\setbox0=\hbox{##2\relax}
\edef\TEST{[\SigText]}
\expandafter\ifx\csname##1\endcsname\TEST\relax\else
\immediate\write16{##1 hat sich geaendert!}\fi
\expandwrite\AUX{\neverexpand\ref{##1}{\TEST}}
\setHI{[\SigText]\ }
\ifnum\HIindent>\maxSpalbr\maxSpalbr\HIindent\fi
\ifnum\HIindent<\altSpalbr\HIindent\altSpalbr\fi
\HI{\marginnote{##1}[\SigText]}
\ifx-##3\relax\else{##3}: \fi
\ifx-##4\relax\else{##4}{\sfcode`.=3000.} \fi
\ifx-##5\relax\else{In: \it ##5}. \fi
\ifx-##6\relax\else{(##6)} \fi\PrII}
\def\PrII##1|Bd:##2|Verlag:##3|Ort:##4|S:##5|J:##6|xxx:##7||{%
\ifx-##1\relax\else{##1} \fi
\ifx-##2\relax\else{\bf ##2}, \fi
\ifx-##4\relax\else{##4:} \fi
\ifx-##3\relax\else{##3} \fi
\ifx-##6\relax\else{##6}\fi
\ifx-##5\relax\else{, ##5}\fi
\ifx-##7\relax\else, \TText{##7}\fi\Par}
\bgroup
\baselineskip12pt
\parskip2.5pt plus 1pt
\hyphenation{Hei-del-berg Sprin-ger}
\sfcode`.=1000
\beginsection References. References

}}

\def\endrefs{%
\expandwrite\AUX{\neverexpand\ref{Spaltenbreite}{\the\maxSpalbr}}
\ifnum\maxSpalbr=\altSpalbr\relax\else
\immediate\write16{Spaltenbreite hat sich geaendert!}\fi
\egroup\write16{Letzte Gleichung: E\the\maxEqNo}
\write16{Letzte Aufzaehlung: I\the\maxItemcount}}

%\L|Abk:|Sig:|Au:|Tit:|Zs:|Bd:|S:|J:|xxx:-||
%\B|Abk:|Sig:|Au:|Tit:|Reihe:|Verlag:|Ort:|J:|xxx:-||
%\Pr|Abk:|Sig:|Au:|Artikel:|Titel:|Hgr:|Reihe:|Bd:|Verlag:|Ort:|S:|J:|xxx:-||

%---Ende makros.mak

\def\v{^\vee}
\def\alq{{\overline\alpha}}
\def\Wf{W_{\mskip-5mu f}}
\def\Hs{\cH^{\rm sph}}
\def\uH{{\underline H}}

%Body

\title{On the Kazhdan-Lusztig basis of a spherical Hecke algebra}

\beginsection Intro. Introduction

Let $\cH$ be an (extended) affine Hecke algebra. It contains the Hecke
algebra $\cH_f$ of the finite Weyl group $\Wf$ as a subalgebra. The
set of elements of $\cH$ which are ``invariant'' under left and right
multiplication by $\cH_f$ is called the spherical Hecke algebra
$\Hs$. The Satake isomorphism identifies $\Hs$ with
$\ZZ[v^{\pm1}][X\v]^{\Wf}$ where $X\v$ is the coweight lattice.

In \cite{KL}, Kazhdan and Lusztig constructed a canonical basis of
$\cH$. This basis is compatible with $\Hs$ and Lusztig has shown,
\cite{Lu2}, that the Kazhdan\_Lusztig elements inside $\Hs$
correspond, under the Satake isomorphism, to the Weyl characters of
the Langlands dual group $G\v$.

The aim of this note is to give a new proof of this result and extend
it to Hecke algebras with unequal parameters. This works as stated, if
the parameters depend only on the root length. If the root system is
of type ${\Ss C}_n$ the situation is more subtle.

The main difference of our approach is that we use Demazure's
character formula while Lusztig used the formula of Weyl. Demazure's
formula is less elementary than Weyl's but apart from that our proof
appears to be simpler. Moreover, Lusztig's proof does not work in the
unequal parameter case since he uses a $q$\_analog of Weyl's
formula, the celebrated Kato-Lusztig formula, which does not seem to
generalize.

\beginsection affineHecke. The extended affine Weyl group

In this and most of the next section, we are setting up notation for
Weyl groups and Hecke algebras and recall the properties which we are
going to use. Proofs can be found, e.g., in Humphreys' book
\cite{Hum}, Macdonald's book \cite{MacBuch} or the nice survey
\cite{NR} of Nelsen and Ram which also presents Lusztig's
approach to our main theorem.

Let $(\Delta_f\subset X,\Delta_f\v\subset X\v)$ be a root datum. The
set of affine roots is $\Delta:=\Delta_f+\ZZ\delta$ whose elements we
regard as affine linear functions on $X\v$ with $\delta$ being the
constant function $1$.  For $\alpha=\alq+m\delta\in\Delta$ let
$$
s_\alpha(\tau):=\tau-\alpha(\tau)\alq\v=
\tau-(\alq(\tau)+m)\alq\v
$$
be the corresponding affine reflection of $X\v$. Let $W^a$ and $\Wf$ be the
groups generated by all reflections $s_\alpha$ with $\alpha\in\Delta$
and $\alpha\in\Delta_f$, respectively. For $\tau\in X\v$ let $t_\tau$
be the translation $t_\tau(\eta):=\eta+\tau$. The group
$W:=\Wf\semidir X\v$ acts on $X\v$ by
$(w,\tau)(\lambda):=wt_\tau(\lambda)=w(\lambda+\tau)$. The group $W^a$
is a subgroup of $W$. More precisely, $W^a=\Wf\semidir Q\v$ where
$Q\v\subseteq X\v$ is the coroot lattice, i.e., the subgroup generated
by $\Delta\v_f$.

We choose a set $\Sigma_f\subseteq\Delta_f$ of simple roots. A root
$\alpha\in\Delta_f$ is called maximal if
$(\alpha+\Sigma_f)\cap\Delta_f=\emptyset$. Clearly, there is one maximal
root for each connected component of the Dynkin diagram. The set
$$
\Sigma:=\Sigma_f\cup\{-\theta+\delta\mid\theta\hbox{ \rm maximal}\}
\subseteq\Delta.
$$
is the set of simple affine roots. The groups $\Wf$ and $W^a$ are
Coxeter groups with generators $\{s_\alq\mid\alq\in\Sigma_f\}$ and
$\{s_\alpha\mid\alpha\in\Sigma\}$, respectively.

Let $\Delta_f^+\subseteq \Delta_f$ and $\Delta^+\subseteq\Delta$ be
the set of positive roots, i.e., those roots which are
$\ZZ_{\ge0}$\_linear combinations of $\Sigma$. With
$\Delta_f^-=-\Delta_f^+$ and $\Delta^-=-\Delta^+$ we have
$\Delta_f=\Delta_f^+\cup\Delta_f^-$, $\Delta=\Delta^+\cup\Delta^-$ and
$$
\Delta^+=(\Delta_f^++\ZZ_{\ge0}\delta)\cup(\Delta_f^-+\ZZ_{>0}\delta).
$$
There is a natural right action of $W$ on the space of affine linear
functions on $X\v$, namely,
$(\alpha^w)(\lambda):=\alpha(w\lambda)$. More precisely, if
$w=t_\tau\wq$ and $\alpha=\alq+m\delta$ then
$$
\alpha^w(\lambda)=\alq^\wq+\alpha(\tau)\delta=\alq^\wq+(\alq(\tau)+m)\delta
$$
Thus, $\Delta$ is stable under $W$ and we define the length of $w\in
W$ as
$$
\ell(w):=\#\{\alpha\in\Delta^+\mid \alpha^w\in\Delta^-\}.
$$
Clearly, $\Omega:=\{w\in W\mid\ell(w)=0\}$ is the stabilizer of
$\Delta^+$ and therefore a subgroup of $W$ and one can show
$W=\Omega\semidir W^a$.  We have $\ell(w^{-1})=\ell(w)$ and
$$
\ell(w)=\|min|\{r\in\NN\mid
\exists\omega\in\Omega,\alpha_1,\ldots,\alpha_r\in \Sigma:
w=\omega s_{\alpha_1}\ldots s_{\alpha_r}\}.
$$
Very useful is the following explicit formula for $\tau\in X\v$,
$w\in \Wf$:
$$6
\ell(w t_\tau)=
\sum_{\alpha\in\Delta_f^+\cap(\Delta_f^+)^w}\big|\alpha(\tau)\big|+
\sum_{\alpha\in\Delta_f^+\setminus(\Delta_f^+)^w}\big|\alpha(\tau)+1\big|.
$$
This implies
$$8
\ell(t_\tau)=\sum_{\alpha\in\Delta_f^+}\big|\alpha(\tau)\big|.
$$
If $\lambda\in X\v_+:=\{\lambda\in
X\v\mid\alpha(\lambda)\ge0\hbox{ for all }\alpha\in\Sigma_f\}$ this
simplifies to
$$2
\ell(t_\lambda)=2\rho(\lambda)\quad\hbox{with}\quad
2\rho:=\sum_{\alpha\in\Delta_f^+}\alpha.
$$

\beginsection HeckeAlgebra. The extended affine Hecke algebra

Let $\cL$ be a ring. For each simple reflection $s=s_\alpha$,
$\alpha\in\Sigma$ we choose an invertible element $v^s\in\cL$ subject
to the condition that $v^{s_1}=v^{s_2}$ if $s_1$, $s_2$ are conjugate
in $W$. In that case we may define $v^w:=v^{s_1}\ldots v^{s_m}$ where
$w=\omega s_1\ldots s_m$ is any reduced expression. We also define
$v^{-w}:=(v^w)^{-1}$. We have mainly two instances of this situation
in mind: first $\cL=\ZZ[v^{\pm s_1},\ldots,v^{\pm s_m}]$ where
$s_1,\ldots,s_m\in\Sigma$ is a set of representatives of $W$\_orbits
in $\Delta$. Secondly, $\cL=\ZZ[v^{\pm1}]$ and $v^s=v^{n_s}$ with
$n_s\in\ZZ$.

Let $\cH$ be the extended Hecke algebra associated to the root datum
$(\Delta_f\subset X,\Delta_f\v\subset X\v)$ and the weights
$v^w$. Thus, $\cH$ is a free $\cL$\_module
with basis $\{H_w\mid w\in W\}$ and relations
$$1
H_wH_{w'}=H_{ww'},\quad\hbox{whenever\ }\ell(ww')=\ell(w)+\ell(w')
$$
and
$$
(H_s-v^{-s})(H_s+v^s)=0,\quad\hbox{for all simple reflections }s=
s_\alpha,\alpha\in\Sigma.
$$
The last relation implies that $H_s$ is invertible and it can be
rephrased as
$$
H_s+v^s=H_s^{-1}+v^{-s}.
$$
The algebra $\cH$ contains the finite dimensional
Hecke algebra $\cH_f:=\oplus_{w\in\Wf}\cL H_w$ attached to the root
system $\Delta_f$ as a subalgebra.

For brevity, we write $H_{t_\tau}=H_\tau$ for $\tau\in X\v$.  Then
equations \cite{E1} and \cite{E2} imply
$$
H_\lambda H_\mu=H_{\lambda+\mu}=H_\mu H_\lambda\quad\hbox{for all}
\quad\lambda,\mu\in X\v_+.
$$
Since every $\tau\in X\v$ is of the form $\lambda-\mu$ with
$\lambda,\mu\in X\v_+$ we can define commuting elements
$$
Y_\tau:=H_\lambda H_\mu^{-1}.
$$
This way, we get a homomorphism
$$
\Phi:\cL[X\v]\pfeil\cH:e^\tau\mapsto Y_\tau.
$$
Moreover, the map
$$12
\cL[X\v]\otimes_\cL\cH_f\Pf\sim\cH:\xi\otimes u\mapsto\Phi(\xi)u
$$
is a linear isomorphism. To see the ring structure on the left-hand
side we have to know the commutation relation between the $H_s$ and
$Y_\tau$. To explain them we need to set up some notation.

\Definition: A simple root $\alpha\in\Sigma_f$ is called {\it special}
if
\Item{} $\alpha$ is the long simple root in a component of
$\Delta_f$ of type ${\Ss C}_n$ ($n\ge1$, with ${\Ss C}_1={\Ss A}_1$)
and
\Item{} $v^{s_{\alpha_0}}\ne v^{s_\alpha}$ where
$\alpha_0:=-\theta+\delta\in\Sigma$ and $\theta\in
W_f\alpha$ is the maximal root.\Par
\noindent In that case we put $v_0^s:=v^{s_{\alpha_0}}$.

\medskip

With this notation we have according to \cite{Lu1}:

\Theorem. Let $\alpha\in\Sigma_f$, $s:=s_\alpha$, and
$\xi\in\cL[X^\vee]$. Then
$$9 H_s\Phi(\xi)-\Phi(s\xi)H_s=
(v^{-s}-v^s)\Phi({\xi-s\xi\over1-e^{-\alpha\v}})
$$
if $\alpha$ is not special and
$$10
H_s\Phi(\xi)-\Phi(s\xi)H_s=
(v^{-s}-v^s)\Phi({\xi-s\xi\over1-e^{-2\alpha\v}})+(v_0^{-s}-v_0^s)
\Phi(e^{-\alpha\v}{\xi-s\xi\over1-e^{-2\alpha\v}})
$$
if it is.

\Remarks: 1. For non\_special $\alpha$ we define
$v_0^s:=v^s$. Then \cite{E9} is a special case of \cite{E10}.\Par
2. If $\alpha\in\Sigma_f$ is special then the
   simple affine root $\alpha_0=-\theta+\delta$ is not $W$\_conjugate to any
   element of $\Sigma_f$. Therefore, the parameter $v^{s_{\alpha_0}}$ is
   possibly different from every parameter $v^{s_\beta}$,
   $\beta\in\Sigma_f$. But it has to occur somewhere in any
   presentation of $\cH$. This explains why not all commutation
   relations can be of the form \cite{E9}.

\medskip

These formulas imply in particular:

\Corollary. The image of $\Phi:\cL[X\v]^{\Wf}\Pfeil\cH$ is in the
center of $\cH$.

\Remark: If the parameters $v^s$ are sufficiently general then one can
show that the image is the entire center, see \cite{Lu1}, but we won't
need this in the sequel.

\medskip
In the definition of $\Phi$ there is nothing special about the
dominant Weyl chamber. For a fixed $w\in\Wf$ we can define
$$
Y_\tau^{(w)}:=H_\lambda H_\mu^{-1}
$$
where $\lambda,\mu\in w(X\v_+)$ with $\tau=\lambda-\mu$. Again, we
get a homomorphism
$$
\Phi_w:\ZZ[X\v]\pfeil\cH:e^\tau\mapsto Y_\tau^{(w)}.
$$
This homomorphism can be expressed in terms of $\Phi$:

\Lemma. For $w\in\Wf$, $\xi\in \cL[X\v]$ holds
$$7
\Phi_w(\xi)=H_w\Phi(w^{-1}\xi)H_w^{-1}.
$$

\Proof: It suffices to prove the formula for $\xi=e^\tau$ where $\tau$
is in the interior of $wX\v_+$. Let $\tau_+=w^{-1}(\tau)\in
X\v_+$. Then formula \cite{E6} implies
$\ell(wt_{\tau_+})=\ell(t_{\tau_+})+\ell(w)$. Hence
$H_{wt_{\tau_+}}=H_wH_{t_{\tau_+}}=H_w\Phi(e^{w^{-1}\tau})$. On the
other hand, formula \cite{E8} implies $\ell(t_\tau)=\ell(t_{\tau_+})$,
and therefore $\ell(t_\tau
w)=\ell(wt_{\tau_+})=\ell(t_\tau)+\ell(w)$. Thus we get
$$
H_w\Phi(e^{w^{-1}\tau})=H_{wt_{\tau_+}}=H_{t_\tau w}
=H_{t_\tau}H_w=\Phi_w(e^\tau)H_w.
$$\qed

\Corollary. The homomorphisms $\Phi$ and $\Phi_w$ coincide on
$\cL[X\v]^{\Wf}$.

Now assume there is an automorphism $x\mapsto\xq$ of $\cL$ with
$\overline{v^w}=v^{-w}$ for all $w\in W$. Then the automorphism
extends uniquely to a duality map $d:\cH\pfeil\cH$ by putting
$$
d(H_w)=H_{w^{-1}}^{-1}\hbox{ for all }w\in W.
$$

\Lemma. Let $w_0$ be the longest element of $\Wf$ and
$\xi\in\ZZ[X\v]$. Then
$$
d(\Phi(\xi))=\Phi_{w_0}(\xi)=H_{w_0}\Phi(w_0\xi)H_{w_0}^{-1}.
$$

\Proof: We may assume $\xi=e^\tau$ with $\tau\in X\v_+$. Then
$$
d(\Phi(e^\tau))=d(Y_\tau)=d(H_{t_\tau})=
H_{t_{-\tau}}^{-1}=(Y_{-\tau}^{(w_0)})^{-1}=Y^{(w_0)}_\tau=\Phi_{w_0}(e^\tau)
$$\qed

\Corollary selfdual. Let $\xi\in\ZZ[X\v]^{\Wf}$ and
$h=\Phi(\xi)$. Then $d(h)=h$.

\beginsection Rightsymmetric. The right spherical submodule

Every right coset in $W/\Wf$ is of the form $t_\tau\Wf$ with a unique
$\tau\in X\v$. It contains a unique element $m_\tau:=t_\tau w_\tau$ of
minimal length. Explicitly, $w_\tau\in \Wf$ is minimal with
$w_\tau^{-1}(\tau)\in -X\v_+$.

The Bruhat order on $W$ induces an order relation on $X\v$:
$$
\lambda\le\mu{\buildrel{\rm def}\over\Longleftrightarrow} m_\lambda\le m_\mu.
$$
\def\ggk{\mathop{\vcenter{\offinterlineskip\hbox{$>$}\vskip2pt\hbox{$=$}\hbox{$<$}}}}
This order relation satisfies (see \cite{Hum}~Prop. 5.7):
$$5
s_\alpha(\tau)\ggk\tau\Longleftrightarrow\alpha(\tau)\ggk0\quad\hbox{for all
}\alpha\in\Delta^+,\tau\in X\v
$$
%$$5
%\hbox{If $\alpha\in\Delta^+$,$\tau\in X\v$ then }
%\cases{
%\alpha(\tau)>0\Rightarrow s_\alpha(\tau)>\tau&\cr
%\alpha(\tau)=0\Rightarrow s_\alpha(\tau)=\tau&\cr
%\alpha(\tau)<0\Rightarrow s_\alpha(\tau)<\tau&\cr}
%$$
and is, in fact, the coarsest order relation with this property.

\Lemma alphastring. For $\tau\in X\v$ and $\alpha\in\Delta_f^+$ with
$N:=\alpha(\tau)\ge0$ let
$$
\tau_0=\tau,\tau_1=\tau-\alpha\v,\ldots,
\tau_N=\tau-N\alpha\v=s_\alpha(\tau)
$$
be the $\alpha\v$\_string through $\tau$. Then
$$
\tau_N>\tau_0>\tau_{N-1}>\tau_1>\tau_{N-2}>\tau_2>
\ldots>\tau_{\lfloor N/2\rfloor}.
$$

\Proof: If $N=0$ there is nothing to show. If $N>0$ we get
$\tau_N>\tau_0$ by \cite{E5}. If $N=1$ we are done, so assume $N>1$
and consider the affine root $\beta=-\alpha+\delta$. Then
$\beta(\tau_{N-1})=N-1>0$ and $s_\beta(\tau_{N-1})=\tau_0>\tau_{N-1}$,
again by \cite{E5}. The remaining inequalities follow by replacing
$\tau$ by $\tau-\alpha\v$.\qed

Consider the following left submodule of $\cH$:
$$
\cM:=\{h\in\cH\mid h H_w=v^{-w}h\quad\hbox{for all }w\in W_f\}.
$$
It is easy to see that $\cM\cap\cH_f=\cL\otheta$ with
$$
\otheta:=\sum_{w\in\Wf}v^{ww_0}H_w\in\cH_f.
$$
Then \cite{E12} implies $\cM\cong\cH\otimes_{\cH_f}\cL\otheta$. Since the
elements $m_\tau$, $\tau\in X\v$, represent the cosets $W/W_f$ we
conclude that the elements
$$
M_\tau:=H_{m_\tau}\otheta=v^{m_\lambda w_0}\sum_{w\in t_\lambda\Wf}
v^{-w}H_w
\quad\hbox{with}\quad\tau\in X\v
$$
form an $\cL$\_basis of $\cM$. On the other hand, the map
$$
\Psi:\cL[X\v]\pfeil\cM:p\mapsto\Phi(p)\otheta
$$
is an isomorphism of $\cL$\_modules. In particular, we obtain a second
basis of $\cM$ namely the elements $\Psi(e^\tau)=Y_\tau\otheta$,
$\tau\in X\v$.

By transport of structure, the Hecke algebra $\cH$ acts also
on $\cL[X\v]$. Explicitly, we have
$$
Y_\eta(e^\tau)=e^{\tau+\eta}\quad\hbox{for }\eta\in X\v
$$
and
$$4
H_s(e^\tau)=v^{-s}e^{s(\tau)}+
(v^{-s}-v^s+(v_0^{-s}-v_0^s)e^{-\alpha\v})
{e^\tau-e^{s(\tau)}\over1-e^{-2\alpha\v}}\quad\hbox{for }
s=s_\alpha, \alpha\in\Sigma_f.
$$

The basis $M_\tau$ of $\cM$ gives rise to a basis
$p_\tau:=\Psi^{-1}(M_\tau)$ of $\cL[X\v]$.

\Lemma demazure. For $\tau\in X\v$ choose any $w\in\Wf$ with
$\tau_+:=w^{-1}(\tau)\in X\v_+$. Then
$$
p_\tau=v^{w_\tau}v^wH_w(e^{\tau_+}).
$$

\Proof: Equation \cite{E7} implies
$$
H_{m_\tau}H_{w_\tau^{-1}}=H_{t_\tau}=\Phi_w(e^\tau)=
H_w\Phi(e^{\tau_+})H_w^{-1}.
$$
Hence
$$
M_\tau=H_w\Phi(e^{\tau_+})H_w^{-1}H_{w_\tau^{-1}}^{-1}\otheta=
v^{w_\tau}v^wH_w\Psi(e^{\tau_+}).
$$\qed

\Lemma triangular. For every $\tau\in X\v$ holds
$p_\tau\in\sum\limits_{\eta\le\tau}\cL e^\eta$.

\Proof: Let $w\in\Wf$ be the shortest element with
$w^{-1}(\tau)\in X\v_+$ and let $w=s_1\ldots s_m$ be a reduced
expression. If $m=0$ then $p_\tau=v^{w_\tau}e^\tau$ and we are
done. For $m\ge1$ put $\tau':=s_1(\tau)$ and $w'=s_1w$. Then we have
$p_\tau=H_{s_1}(p_{\tau'})$. By induction we may assume that every
monomial $e^\eta$ occurring in $p_{\tau'}$ satisfies $\eta\le\tau'$. The
monomials $e^{\eta'}$ occurring in $H_{s_1}(e^\eta)$ are all in the
$\alpha_1$\_string with endpoint $\eta$. If $\alpha(\eta)\le0$ then
\cite{alphastring} implies $\eta'\le\eta\le\tau'<\tau$ and we are
done. If $\alpha(\eta)>0$ the same holds except for
$\eta'=s_1(\eta)>\eta$. But then $\tau=s_1(\tau')>\tau'\ge\eta$
implies $\tau\ge\eta'$.\qed

To cover special simple reflections, we introduce the root system
$\tilde\Delta_f\subseteq X$ which is generated by
$\tilde\Sigma_f:=\{\epsilon(\alpha)\alpha\mid\alpha\in\Sigma_f\}$ with
$$
\epsilon(\alpha):=\cases{
{1\over2}&if $\alpha$ is special;\cr
1&otherwise.\cr}
$$
Correspondingly, $\tilde\Delta\v_f$ is generated by
$\tilde\Sigma\v_f:=\{\epsilon(\alpha)^{-1}\alpha\v\mid\alpha\v\in\Sigma\v_f\}$.
In other words, $\tilde\Delta\v_f=\Delta\v_f$ if none of the simple
roots are special while $\tilde\Delta\v_f={\Ss C}_n$ if
$\Delta_f={\Ss C}_n$ and the long root is special.

We recall the Demazure operators (see, e.g., \cite{Dem}). For each
simple reflection $s=s_\alpha$, $\alpha\in\tilde\Sigma_f$ we define
$$
\Delta_s:=s+(1-e^{-\alpha\v})^{-1}(1-s)
$$
which acts on $\ZZ[X\v]$. If $w=s_1\ldots s_m\in\Wf$ is a reduced
expression then $\Delta_w=\Delta_{s_1}\ldots\Delta_{s_m}$ depends only
on $w$. For $w\in\Wf$ and $\lambda\in X\v_+$ the element
$\Delta_w(e^\lambda)$ is called a Demazure character. We parameterize
it as follows: for $\tau\in X\v$ let $w\in\Wf$ be such that
$\tau_+=w^{-1}(\tau)\in X\v_+$. Then put
$\delta_\tau:=\Delta_w(e^{\tau_+})$. This does not depend on the choice
of $w$.

Now we can be more specific about the coefficients in
\cite{triangular}.

\Lemma integral. Let $\cL_{++}\subseteq\cL$ be a non\_unital subring
which contains all $v^s$ where $s:=s_\alpha$, $\alpha\in\Sigma_f$
and, moreover, $v^sv_0^{\pm s}$ in case $\alpha$ is special.  Let
$\tau\in X\v$. Then
$$
p_\tau\in v^{w_\tau}\big(\delta_\tau+\cL_{++}[X\v]\big).
$$

\Proof: Let again $w\in\Wf$ with $\tau_+=w^{-1}(\tau)\in
X\v_+$. Choose a reduced expression $w=s_1\ldots s_m$. Then, by
\cite{demazure},
$$
p_\tau=v^{w_\tau}(v^{s_1}H_{s_1})\ldots(v^{s_m}H_{s_m})(e^{\tau_+}).
$$
By \cite{E4} the operator $v^sH_s$, $s=s_\alpha, \alpha\in\Sigma_f$
can be expressed as
$$
v^sH_s=\Delta_s-(v^s)^2(1-e^{-\alpha\v})^{-1}(1-s)
$$
if $s$ is not special and
$$
v^sH_s=\Delta_s-[(v^s)^2-(v^sv_0^{-s}-v^sv_0^s)e^{-\alpha\v}]
(1-e^{-2\alpha\v})^{-1}(1-s)
$$
if $s$ is special. This implies the assertion.\qed

\beginsection Spherical. The spherical Hecke algebra

Every double coset in $\Wf\backslash W/\Wf$ is of the form $\Wf
t_\lambda \Wf$ with a unique $\lambda\in X\v_+$. It contains a unique
longest element namely $n_\lambda:=w_0t_\lambda$. Thus we have
$n_\lambda=m_{w_0\lambda}w_0$.

We put
$$
\eqalign{\Hs:
&=\{h\in\cM\mid H_wh=v^{-w}h\quad\hbox{for all }w\in W_f\}=\cr
&=\{h\in\cH\mid H_wh=hH_w=v^{-w}h\quad\hbox{for all }w\in W_f\}.\cr}
$$
Since $\cM=\cH\otheta$, we obtain $\Hs=\otheta\cH\cap\cH\otheta$. The
bijection $\Psi:\cL[X\v]\pf\sim\cM$ induces a bijection
$$11
\Psi:\cL[X\v]^{\Wf}\pf\sim\Hs:p\mapsto\Phi(p)\otheta=\otheta\Phi(p),
$$
called the {\it Satake isomorphism}.

\Remark: Let $P:=v^{-w_0}\sum_{w\in\Wf}(v^{w})^2\in\cL$.
Then $\otheta^2=P\otheta$ which implies
$\Psi(\xi_1\xi_2)=P\Psi(\xi_1)\Psi(\xi_2)$. Hence, if $P$ is not a
zero divisor, we can define a new multiplication $h_1\ast h_2:={1\over
P}h_1h_2$ on $\Hs$ for which $\Hs$ becomes a commutative ring with
identity element $\otheta$ and \cite{E11} is an isomorphism of rings.

\medskip
For every $\lambda\in X\v_+$ put
$$3
N_\lambda:=\sum_{\tau\in\Wf\lambda}v^{w_\tau}M_\tau=
v^{n_\tau}\sum_{w\in\Wf t_\lambda\Wf}v^{-w}H_w.
$$
Then the
$N_\lambda$ form an $\cL$\_basis of $\cH^{\rm sph}$. Via $\Psi$, they
give rise to a basis $P_\lambda:=\Psi^{-1}(N_\lambda)$ of
$\cL[X\v]^{\Wf}$. For the root system ${\Ss A}_{n-1}$ they are
basically the Hall\_Littlewood polynomials.

For $\lambda\in X\v_+$ we denote the Demazure character
$\delta_{w_0(\lambda)}=\Delta_{w_0}(e^\lambda)$ by $s_\lambda$. It is
well known that the $s_\lambda$ form a $\ZZ$\_basis of
$\ZZ[X\v]^{\Wf}$. For the root system ${\Ss A}_{n-1}$ they are
basically the Schur polynomials.

\Theorem sphertriang. Let $\cL_{++}\subseteq\cL$ be as in
\cite{integral} and let $\lambda\in X\v_+$. Then
$$
P_\lambda\in s_\lambda+\sum_{\mu\in X\v_+\atop\mu<\lambda}\cL_{++}s_\mu.
$$

\Proof: The definition \cite{E3} and \cite{integral} imply
$P_\lambda\in s_\lambda +r_\lambda$ with
$r_\lambda\in\cL_{++}[X\v]$. Since $s_\lambda$ is $\Wf$\_invariant we
also have $r_\lambda\in\cL_{++}[X\v]^{\Wf}$. Thus $r_\lambda$ is a
$\cL_{++}$\_linear combination of $s_\mu$'s. Finally,
\cite{triangular} implies that every $s_\mu$ occurring in $r_\lambda$ has
$\mu<\lambda$.\qed

\beginsection KazhdanLusztig. Kazhdan-Lusztig elements

In \cite{KL}, Kazhdan and Lusztig constructed their celebrated basis of
$\cH$. Recall, that $\cL$ was supposed to be equipped with an
involution $x\pfeil\xq$ such that $\overline{v^w}=v^{-w}$ for all
$w\in W$. Moreover, fix an additive subgroup $\cL_{++}\subseteq\cL$ and
put $\cH_{++}:=\sum_{w\in W}\cL_{++}H_w\subseteq\cH$.

\Definition: A {\it KL\_element} for $w\in W$ is an element
$\uH_w\in\cH$ with
\Item{}$d(\uH_w)=\uH_w$ and
\Item{}$\uH_w\in H_w+\cH_{++}$.

\medskip\noindent As for existence and uniqueness, we have the
following theorem. Its proof is quite easy and can be found in
\cite{LuBuch}. See also a revised version of \cite{Soe} on the
Soergel's homepage.

\Theorem. For $\cL^{\pm}:=\{x\in\cL\mid\xq=\pm x\}$ consider the homomorphism
$\phi:\cL_{++}\pfeil\cL^-:x\mapsto x-\xq$.
\Item{} Assume that $\phi$ is injective, i.e., $\cL_{++}\cap\cL^+=0$. Then
for every $w\in W$ there is at most one KL\_element $\uH_w$.
\Item{} Assume that $\phi$ is surjective. Then for every $w\in W$ there is
a KL\_element $\uH_w$ which is even triangular, i.e.,
$\uH_w\in\sum_{v\le w}\cL H_v$.\Par

We come to the main theorem of our paper where we explicitly compute
the KL\_elements for $n_\lambda$. This generalizes a result of
Lusztig \cite{Lu2} who proved the same in case of equal
parameters. As mentioned in the introduction, his proof is quite
different from ours.

\Theorem. Let $\cL_{++}\subseteq\cL$ be a non\_unital
subring which contains all $v^s$ where $s:=s_\alpha$, $\alpha\in\Sigma_f$
and moreover $v^sv_0^{\pm s}$ in case $\alpha$ is special. Let $\lambda\in
X\v_+$. Then $\Psi(s_\lambda)$ is a KL\_element for $n_\lambda$.

\Proof: We verify that $\Psi(s_\lambda)$ satisfies the defining
properties of $\uH_{n_\lambda}$.

First, we have $\cM\cap\cH_f=\cL\otheta$. Hence
$d(\otheta)\in\cL\otheta$ and therefore $d(\otheta)=\otheta$. Together
with \cite{selfdual} this implies that all elements of
$\Psi(\ZZ[X\v]^{\Wf})$, in particular $\Psi(s_\lambda)$, are selfdual.

The spherical algebra $\Hs$ has two bases, namely $N_\lambda$ and
$\Psi(s_\lambda)$ with $\lambda\in X\v_+$. By \cite{sphertriang}, the
transition matrix from the former to the latter is unitriangular with
nondiagonal coefficients in $\cL_{++}$. Thus the same holds for its
inverse. This shows
$$
\Psi(s_\lambda)\in N_\lambda+\sum_\mu \cL_{++}N_\mu\subseteq
H_{n_\lambda}+\cH_{++}.
$$\qed

\noindent{\bf Remark:} The most important case is the one considered
by Lusztig: here $\cL=\ZZ[v,v^{-1}]$ with $\vq=v^{-1}$ and
$\cL_{++}=v\ZZ[v]$. Then $\phi:\cL_{++}\pfeil\cL^-$ is bijective which
implies that all KL\_elements exist and are unique. Moreover, the
parameters are of the form $v^s=v^{n_s}$ with $n_s\in\ZZ$. Then the
conditions of the Theorem boil down to $n_s>0$ for $s=s_\alpha$,
$\alpha\in\Sigma_f$ and $|n_{s_0}|<n_s$ if $\alpha$ is special and
$s_0$ is the associated affine reflection. In particular, $n_{s_0}$
may be negative.

\Corollary. Assume additionally, that $\phi:\cL_{++}\pfeil\cL^-$ is
injective. For $\lambda\in X\v_+$ let $L_\lambda$ be the irreducible
$\tilde G\v$\_module with highest weight $\lambda$. Let $m_{\lambda\mu}^\nu$
be the multiplicity of $L_\nu$ in $L_\lambda\otimes L_\mu$. Then
$$
\uH_{n_\lambda}\ast\uH_{n_\mu}
=\sum_{\nu\in X\v_+}m_{\lambda\mu}^\nu\uH_{n_\nu}.
$$

\Proof: This expresses the fact that $s_\lambda$ is the character
of $L_\lambda$ (Demazure's character formula).\qed

\noindent Assume we have a (evaluation) homomorphism $\epsilon:\cL\pfeil\ZZ$
with $\epsilon(v^w)=1$ for all $w\in W$. Then

\Corollary. Let $\uH_w=\sum_{u\in W}p_{uw}H_u$. For $\lambda\in X,\mu\in
X\v_+$ let $L_\lambda(\mu)$ be the $\mu$\_weight space in the
irreducible representation of $\tilde G\v$ with highest weight
$\lambda$. Then $\|dim|L_\lambda(\mu)=\epsilon(p_{n_\mu n_\lambda})$.

\Proof: This uses the fact that for $v^s=1$, the Hecke algebra
degenerates to the group algebra of $W=X\v\dirsemi W_f$ and $\Psi$ becomes
the ``obvious'' map $e^\tau\mapsto\sum_{w\in W_f}t_\tau w$.\qed

\beginrefs

\L|Abk:Dem|Sig:De|Au:Demazure, Michel|Tit:D\'esingularisation des
vari\'et\'es de Schubert g\'en\'eralis\'ees|Zs:Ann. Sci. \'Ecole Norm.
Sup.|Bd:7|S:53--88|J:1974|xxx:-||

\B|Abk:Hum|Sig:Hum|Au:Humphreys, J.|Tit:Reflection groups and Coxeter
groups|Reihe:Cambridge Studies in Advanced Mathematics,
{\bf29}|Verlag:Cambridge University Press|Ort:Cambridge|J:1990|xxx:-||

\L|Abk:KL|Sig:KL|Au:Kazhdan, D.; Lusztig, G|Tit:Representations of Coxeter
groups and Hecke
algebras|Zs:Invent. Math.|Bd:53|S:165--184|J:1979|xxx:-||

\Pr|Abk:Lu2|Sig:\\Lu|Au:Lusztig, G.|Artikel:Singularities, character
formulas, and a $q$-analog of weight multiplicities|Titel:Analysis and
topology on singular spaces, II, III (Luminy, 1981)|Hgr:-|Reihe:Ast\'erisque%
|Bd:101-102|Verlag:Soc. Math. France|Ort:Paris|S:208--229|J:1983|xxx:-||

\L|Abk:Lu1|Sig:\\Lu|Au:Lusztig, G.|Tit:Affine Hecke algebras and their
graded version|Zs:J. Amer. Math. Soc.|Bd:2|S:599--635|J:1989|xxx:-||

\B|Abk:LuBuch|Sig:\\Lu|Au:Lusztig, G.|Tit:Introduction to quantum
groups|Reihe:Progress in Mathematics {\bf
110}|Verlag:Birkh\"auser|Ort:Boston|J:1993|xxx:-||

\B|Abk:MacBuch|Sig:Mac|Au:Macdonald, I.|Tit:Affine Hecke algebras and
orthogonal polynomials|Reihe:Cambridge Tracts in
Mathematics|Verlag:Cambridge University Press|Ort:Cambridge|J:2003|xxx:-||

\Pr|Abk:NR|Sig:NR|Au:Nelsen, K.; Ram, A.|Artikel:Kostka-Foulkes
polynomials and Macdonald spherical functions|Titel:Surveys in
Combinatorics 2003|Hgr:C.D. Wensley ed.|Reihe:London
Math. Soc. Lecture Note Ser.|Bd:307|Verlag:Cambridge
Univ. Press|Ort:Cambridge|S:325--370|J:2003|xxx:math.RT/0401298||

\L|Abk:Soe|Sig:Soe|Au:Soergel, W.|Tit:Kazhdan-Lusztig-Polynome und eine
Kombinatorik f\"ur Kipp-Moduln|Zs:Represent.
Theory|Bd:1|S:37--68|J:1997|xxx:-||

\endrefs

\bye